\documentclass[amstex]{article}
\usepackage{amssymb}
\usepackage{latexsym}
\usepackage{times}

\newcommand{\ds}{\displaystyle}

\newcommand{\cbi}[2]{\left(\begin{array}{c}#2\\#1\end{array}\right)}

\newtheorem{dfn}{Definition}[section] 
\newtheorem{rmk}{Remark}[section]
 
\newtheorem{cor}{Corollary}[section]
\newtheorem{prop}{Proposition}[section] 
\newtheorem{lem}{Lemma}[section]

\newcommand{\Pf}{{\em Proof}. }
\newcommand{\EPf}
{%
\mbox{}%
\nolinebreak%
\hfill%
\rule{2mm}{2mm}%
\medbreak%
\par%
}

\newcommand{\R}{\mathbb R}

\newcommand{\U}{\mathcal{U}{}}

\newcommand{\w}{\mathcal{W}{}}

\newcommand{\s}{\mathcal{S}{}} 

\renewcommand{\d}{\mathcal{D}{}} 
\newcommand{\h}{\mathcal{H}{}}

\newcommand{\ad}{\mbox{ad}}

\newcommand{\rad}{\mbox{rad}}

\newcommand{\stwovd}{\mbox{$S^2(V^{\star})$}}
\newcommand{\srvd}{\mbox{$S^r(V^{\star})$}}

\newcommand{\Sh}{\mbox{$\mathcal{S}(\mathcal{H})$}}
\newcommand{\Wh}{\mbox{$\mathcal{W}(\mathcal{H})$}}

\newcommand{\ZWh}{\mbox{$\mathcal{Z}\mathcal{W}(\mathcal{H})$}}
\newcommand{\ZW}{\mbox{$\mathcal{Z}\mathcal{W}$}}

\newcommand{\Soneh}{\mbox{$\mathcal{S}^{(1)}(\mathcal{H})$}}

\newcommand{\Srh}{\mbox{$\mathcal{S}^{(r)}(\mathcal{H})$}}
\newcommand{\OM}{\mbox{$\Omega$}}

\newcommand{\D}{\nabla}

\begin{document}
\title{On the geometry of the characteristic class of a star product on a symplectic 
manifold}
\author{Pierre Bieliavsky
and Philippe Bonneau\footnote{Research supported by 
the Communaut\'e fran\c caise de Belgique, through an
Action de Recherche Concert\'ee de la Direction de la Recherche 
Scientifique}}
\date{\today}

 \maketitle
\begin{abstract}
The characteristic class of a star product on a symplectic manifold
appears as the class of a deformation of a given symplectic 
connection, as 
described by Fedosov. In contrast, one usually 
thinks of the characteristic class of a star product as the class of a 
deformation of the Poisson 
structure (as in Kontsevich's work).
In this paper, we present, in the symplectic framework, a natural
procedure for constructing a star product 
by directly quantizing a deformation of the symplectic structure. 
Basically, in Fedosov's recursive formula for the star product with {\bf zero} characteristic class,
we replace the symplectic structure by one of its formal deformations in the parameter $\hbar$.
We then show that every equivalence class of star products contains such an element. 
Moreover, within a given class, equivalences between such star products are realized 
by formal one-parameter families of diffeomorphisms, as produced by Moser's argument.
\end{abstract}

\section{Introduction}
The characteristic class of a star product on a symplectic manifold
appears as the class of a deformation of a given symplectic 
connection, as 
described by Fedosov \cite{Fedbook,Fedb94a}. In contrast, one usually 
thinks of the characteristic class of a star product as the class of a 
deformation of the Poisson 
structure  \cite{Konm97b}.
In this paper, we present, in the symplectic framework, a natural
procedure for constructing a star product 
by directly quantizing a deformation of the symplectic structure. 
Basically, in Fedosov's recursive formula for the star product with {\bf zero} characteristic class,
we replace the symplectic structure by one of its formal deformations in the parameter $\hbar$.
We then show that every equivalence class of star products contains such an element. 
Moreover, within a given class, equivalences between such star products are realized 
by formal one-parameter families of diffeomorphisms, as produced by Moser's argument.
More precisely, let $(M,\omega)$ be a compact symplectic manifold. Let 
$\{ \Omega{(t)}\}_{t\in]-\epsilon,\epsilon[}$
be a smooth path of symplectic structures on $M$ such that $\Omega_{0}=\omega$. 
The pair $(M,\{\Omega{(t)}\})$  defines a regular Poisson structure $\hat{\Omega}$ on 
$\hat{M}=M\times
]-\epsilon,\epsilon[$ whose symplectic leaves are $\{(M\times\{t\},\Omega{(t)})\}$. 
Applying Fedosov's method to 
$(\hat{M},\hat{\Omega})$, one obtains a tangential star product $\hat{\star}$ on 
$(\hat{M},\hat{\Omega})$ with 
zero characteristic class. The ``infinite jet at $0$ of $\hat{\star}$ in $t=\hbar$'' 
then defines a star product $\star$ on $(M,\omega)$ to which is 
associated the de Rham class 
$[\Omega^\hbar]_{\mbox{de Rham}}$ 
where $\Omega^\hbar$ denotes the infinite jet  at $0$ of $\{\Omega{(t)}\}$ in $t=\hbar$. 
If $\{\Omega'{(t)}\}$ is 
such that 
$[\Omega'{(t)}]_{\mbox{de Rham}}=[\Omega{(t)}]_{\mbox{de Rham}}\,\forall t$, 
then an equivalence between the 
corresponding star products is realized as the infinite jet 
of a family of diffeomorphisms $\{\varphi_t\}$ ---whose existence is 
guaranteed by Moser's argument--- such that $\varphi_t^\star\Omega'{(t)}=\Omega{(t)}$.

This work is motivated by the question of obtaining a quantum analogue of Kirwan's map
when considering the problem of commutation between Marsden-Weinstein reduction and 
deformation quantization. However this point is not investigated in the 
present article.

\section{Fedosov construction on regular Poisson manifolds}\label{FEDO}
We present Fedosov star products on regular Poisson manifolds 
\cite{Fedbook,Fedb94a}
by mean of a partial connection defined (only) on the characteristic 
distribution of the Poisson structure. By this we avoid considering 
Poisson affine connections (cf. Lemma \ref{CONN}). This little point excepted, there is 
essentially nothing new in the present section. But it sets the notations 
and presents Fedosov's construction in a completely intrinsic way.
\subsection{Linear Weyl algebra}
Let $(V,\omega)$ be a real symplectic vector space and consider the 
associated Heisenberg Lie algebra $\h$ over the dual space $V^\star$. That is 
$\h=V^\star\oplus\R\hbar$ where $\hbar$ is central and where the Lie 
bracket of two elements $y,y'\in V^\star$ is defined by 
$[y,y']=y'({}^\sharp y)\hbar$. The map $V^\star\stackrel{\sharp}{\to} V$ 
being the isomorphism induced by $\omega$. Denote  by $\s(\h)$ (resp. $\U(\h)$) 
the symmetric (resp. the universal enveloping) algebra of $\h$ and 
consider the complete symmetrization map 
$\s(\h)\stackrel{\varphi}{\to}\U(\h)$ given by the 
Poincar\'e-Birkhoff-Witt 
theorem. The symmetric product on $\s(\h)$ will be denoted by 
$\bullet$, 
while $\star$ will denote the product on $\s(\h)$ transported via $\varphi$ of 
the universal product on $\U(\h)$. 
\begin{lem}\label{GRAD}
There exists one and only one grading $\Sh =: \oplus_{r \geq 0} \Srh$ on 
$\Sh$ such that~:
$$
\begin{array}{cccc}
(i) & \srvd & \subset & \Srh \\
(ii) & \Srh \star \mathcal{S}^{(s)}(\mathcal{H}) & \subset & \mathcal{
S}^{(r+s)}(\mathcal{H}),
\end{array}
$$
where $\srvd$ denotes the $r$-th  symmetric power of $V^{\star}$.
This grading is compatible with the symmetric product $\bullet$ as well.
\end{lem}
One then defines the {\bf linear Weyl algebra} $\Wh$ as the direct product 
$\Wh := 
\prod_{r=0}^{\infty} \Srh$ endowed with the extended product $\star$. Note 
that the symmetric product $\bullet$ extends to $\Wh$ as well. The center 
$\ZWh$ of $(\Wh,\star)$ is canonically isomorphic to the space of power series 
$\R[[\hbar]]$.

By using the symplectic structure, one gets an identification between the 
Lie algebra $sp(V,\omega)$ and the second symmetric power $\stwovd$~:
$$
\begin{array}{ccc}
sp(V,\omega)&\longrightarrow&\stwovd\subset\Wh\\
A&\mapsto&\underline{A}.
\end{array}
$$
\begin{lem}\label{HAM} For all $a\in \Wh$ and $A\in sp(V,\omega)$, one has
$$[\underline{A},a]=2\hbar \, A(a),$$
where $[\, , \,]$ denotes the Lie bracket on $\Wh$ induced by the 
associative product $\star$.
\end{lem}
\Pf Both $\ad(\underline{A})$ and $\hbar A$ are derivations 
of $(\Wh,\star)$. It is therefore sufficient to verify formula (i) on 
generators.
\EPf

\noindent The isomorphism $V\stackrel{\flat}{\to} V^{\star}$ defines an injection
$V  \stackrel{\mu}{\to}  \Wh $ which we call the {\bf linear moment}.
Observe that, viewed as an element of $\Wh\otimes V^\star$, $\mu$ is fixed 
under the action of the symplectic group $Sp(V,\omega)$.

\noindent Both products $\star$ and $\bullet$ extend naturally to the space 
$\Wh\otimes\Lambda^\bullet(V^\star)$ of multilinear forms on $V$ valued in $\Wh$.
We define the total degree $t$ 
of an element $a\otimes \omega , \, a \in \Srh, \omega \in 
\Lambda^p(V^\star)$ by $t=p+r$. With respect to this degree on 
$\Wh\otimes\Lambda^\bullet(V^\star)$, the extended multiplications, again 
denoted by $\star$ and $\bullet$, are graded. The bracket $[\, , \,]$ 
mentioned in Lemma \ref{HAM} therefore extends to $\Wh\otimes\Lambda^\bullet(V^\star)$ 
as well, and, $(\Wh\otimes\Lambda^\bullet(V^\star), [\, , \,])$ is a 
graded Lie algebra.

\noindent To an element $a\otimes x \in \Wh \otimes  \Lambda^p(V)$, one can 
associate the operator 
$$
i_{a\otimes x}:\Wh\otimes\Lambda^\bullet(V^\star) \to \Wh\otimes\Lambda^{\bullet-p}(V^\star),
$$
defined by 
$$
i_{a\otimes x}(b\otimes \omega):=a\bullet b \otimes i_x\omega,
$$
where $i_x\omega$ denotes the usual interior product.
Using the universal property, one gets a map 
$$
\begin{array}{ccc}
(\Wh\otimes V)\times \Wh\otimes\Lambda^\bullet(V^\star) & 
\to & \Wh\otimes\Lambda^{\bullet-p}(V^\star) \\
(X,s) & \mapsto & i_Xs.
\end{array}
$$
In the case where $p$ is odd, since $i_X$ acts ``symmetrically" on the 
``Weyl part" and ``an\-ti-sym\-met\-rical\-ly" on the ``form part", one 
has $i_X^2=0$.\\
In the same way, if $Y \subset \Wh$ is a subspace such that $[Y,Y]
\subset \mathcal{Z}\Wh$ (e.g. $Y=\Soneh$), to any element $U\in 
Y\otimes\Lambda^p(V^\star)$, one can associate the operator $$
\ad(U) : \Wh\otimes\Lambda^\bullet(V^\star)  \to  \Wh\otimes\Lambda^{\bullet+p}(V^\star).
$$
Using  Jacobi identity on the ``Weyl part", one observes that, if $p$ is 
odd, 
one has $\ad(U)^2=0$.
\begin{dfn}
Using the duality 
$$
\begin{array}{ccc}
\Soneh \otimes V^\star & \to & \Soneh \otimes V \\
U &\mapsto & {}^\sharp U,
\end{array}
$$
one defines the cohomology (resp. homology) operator $\delta$ (resp. 
$\delta^\star$) by 
$$
\begin{array}{ccc}
\hbar\delta & := & \ad(\mu)\\
\delta^\star &:= & i_{{}^\sharp\mu},
\end{array}
$$
where the linear moment $\mu$ is viewed as an element of $\Soneh \otimes 
V^\star$. 
For a form $a \in \Wh\otimes\Lambda^\bullet(V^\star)$ with total degree $t$, we set 
$$
\delta^{-1}a:=\left[ 
                        \begin{array}{ccc}
                        \frac{1}{t}\delta^\star a & if & t>0 \\
                        0 & if & t=0.
                        \end{array}
\right.
$$
One extends this definition to the whole $\Wh\otimes\Lambda^\bullet(V^\star)$.
\end{dfn}
\begin{lem}\label{HODGE}(``Hodge decomposition")
$$
\delta\delta^{-1}+\delta^{-1}\delta=Id-\mbox{pr}_{0} 
$$
where $\mbox{pr}_{0}$ is the canonical projection $\Wh\otimes\Lambda^\bullet(V^\star) 
\stackrel{\mbox{pr}_{0}}{\to} \mathcal{Z}\Wh$.
\end{lem}
\Pf
We  observe that $\delta$ and $\delta^\star$ are anti-derivations of 
degree $+1$ of $(\Wh\otimes\Lambda^\bullet(V^\star),\, {\bullet} \,)$. Their 
anti-commutator being a derivation of degree $0$, it is therefore 
sufficient to check the formula on generators.
\EPf
\noindent Observe that $\delta$ is an anti-derivation of degree $+1$ of 
$(\Wh\otimes\Lambda^\bullet(V^\star),\, {\star} \,)$.

\subsection{The Weyl bundle}
Let $(N,\Lambda)$ be a regular Poisson manifold. The Poisson bivector 
$\Lambda$ induces a short sequence of vector bundles over $N$~:
$$
0\to\rad(\Lambda)\to T^\star(N) \stackrel{\iota^\star}{\to}\d^\star\to 0
$$
where $\d\stackrel{\iota}{\to}T(N)$ denotes the characteristic distribution 
associated to $\Lambda$ \cite{Vais94a}, and where $\rad(\Lambda)$ is the radical 
of $\Lambda$ in $T^\star(N)$. One therefore gets a non-degenerate foliated 
2-form $\omega^\d\in\Omega^2(\d)$, dual to the canonical one on the quotient
$T^\star(N)/\rad(\Lambda)=\d^\star$. Fix a rank($\d$)-dimensional 
symplectic vector space $(V,\omega)$, and, for all $x\in N$, define
$$
P_x=\{b\in\mbox{Hom}_\R(V,\d_x)\, |\, b^\star\omega^\d_x=\omega\}.
$$
Then $P=\cup_{x\in N}P_x$ is naturally endowed with a structure of 
$Sp(V,\omega)$-principal bundle over $N$ (analogous to the symplectic 
frames in the symplectic case, except that here, one does not have a 
$G$-structure in general).
\begin{dfn}
The {\bf Weyl bundle} is the associated bundle 
$$
\w=P\times_{Sp(V,\omega)}\Wh,
$$
where $\Wh$ is the vector space underlying the linear Weyl algebra defined 
from the data of $(V,\omega)$.
\end{dfn}
The space of p-forms with values in the sections of $\w$ is denoted by 
$\OM^p(\w)$; it is canonically isomorphic to the space of sections of 
the associated bundle 
$P\times_{Sp(V,\omega)}(\Wh\otimes\Lambda^p(V^\star))$. The 
$Sp(V,\omega)$-invariance, at the linear level, of both product $\star$ and 
$\bullet$ on $\Wh\otimes\Lambda^\bullet(V^\star)$ provides graded 
products, again denoted by $\star$ and 
$\bullet$, on $\OM^\bullet(\w)$. In the same way, the operators $\delta$ and 
$\delta^{-1}$ on $\Wh\otimes\Lambda^\bullet(V^\star)$ define operators on 
sections~:
$$
\begin{array}{ccc}
\OM^\bullet(\w) & 
        \begin{array}{c}
        \delta\\
        \longrightarrow\\
        \longleftarrow\\
        \delta^{-1}
        \end{array} &
\OM^{\bullet+1}(\w),    
\end{array}
$$
leading to a Hodge decomposition of sections as in Lemma \ref{HODGE}.
Notes that the bundle $\ZW=P\times_{Sp(V,\omega)}\ZWh$ being 
trivial, its space of sections is isomorphic to $\mathcal{C}^\infty(N)[[\hbar]]$.
\begin{rmk}\label{WB}
{\rm 
Observe that, as a vector bundle, $\w$ is defined as soon as the distribution 
$\d$ is given (cf. Lemma~\ref{GRAD}). The full data of the Poisson
tensor $\Lambda$  
is only needed to define the algebra structure on its space of sections.
}
\end{rmk}
\subsection{Fedosov-Moyal star products}
\begin{dfn}
A {\bf foliated connection} is a linear map 
\begin{displaymath}
\begin{array}{llll}
\nabla : & \underline{\mathcal{D}}\otimes \underline{\mathcal{D}} 
& \longrightarrow & \underline{\mathcal{D}} \\
& u\otimes v & \longmapsto & \nabla_uv
\end{array}
\end{displaymath}
verifying ($f\in \mathcal{C}^\infty (N)$)
\begin{enumerate}
\item[(i)] $\nabla_{fu}v=f\nabla_uv$,
\item[(ii)] $\nabla_ufv=f\nabla_uv + L_{\iota(u)}f \ v$.
\end{enumerate}

\noindent A foliated connection is said to be {\bf symplectic} if 
\begin{enumerate}
\item[(iii)] $\nabla_uv - \nabla_vu - [u,v]=0$,
\item[(iv)] $ \nabla \omega = 0$.
\end{enumerate}
\end{dfn}

\begin{lem}\label{CONN}
On a regular Poisson manifold, a symplectic foliated connection always 
exists.
\end{lem}
\Pf
Choose any linear connection $\D^{0}$ in the vector bundle $\d\to N$. 
Since $\d$ is an involutive tangent distribution, the torsion $T^0$ of the 
connection is well defined as a section of $\d^\star\otimes\mbox{End}(\d)$. 
One then obtains a ``torsion-free" connection $\D^1=\D^0-\frac{1}{2}T^0$ 
in $\d$. Now, the formula 
$$
\omega^\d(S(u,v),w)=\frac{1}{3}\left( 
\D^1_u\omega^\d(v,w)+\D^1_v\omega^\d(u,w)
\right)
$$
defines a tensor $S$, section of $\d^\star\otimes\d^\star\otimes\d$ such 
that 
$\D=\D^1+S$ is as desired.
\EPf

\noindent Now, fix such a foliated symplectic connection $\D$ in $\d$ and consider 
its associated covariant exterior derivative 
$$
\OM^p(\w)\stackrel{\partial}{\longrightarrow}\OM^{p+1}(\w),
$$
defined by 
$$
\partial s (u_1,..., u_{p+1}) = \sum_{i=1}^{p+1} (-1)^{i-1} (\nabla_{u_i} 
s) (u_1,..., \hat{u_i},..., u_{p+1}).
$$
Lemma \ref{HAM} then provides a 2-form 
$\underline{R}\in\OM^2(\d\otimes\d)\subset\OM^2(\w)$ defined by the formula
$$
2\hbar\partial^2=\ad(\underline{R}).
$$
\noindent Inductively on the degree, one sees (\cite{Fedbook} (Theorem 5.2.2)) that the equation
$$
\underline{R} +2\hbar(\partial\gamma-\delta\gamma+\gamma^2)=0.
$$
has a unique solution $\gamma\in\OM^1(\w)$ such that 
$\delta^{-1}\gamma=0$.

This implies that the graded derivation
$$
D=\partial-\delta+\ad(\gamma)
$$
of $(\OM^\bullet(\w), \star)$ is flat i.e. $D^2=0$. One then proves, again inductively, that
the projection
$$
\w_D\stackrel{\mbox{pr}_0}{\longrightarrow}\mathcal{C}^\infty(N)[[\hbar]],
$$
where $\w_D$ is the kernel of $D$ restricted to the sections of $\w$, is a linear isomorphism. 
The space of flat sections $\w_D$ being a subalgebra of the sections of 
$\w$ with respect to the product $\star$ ($D$ is a derivation), the above 
linear isomorphism yields a star product on $\mathcal{C}^\infty(N)$ called 
{\bf Fedosov-Moyal} star product on $(N,\Lambda)$.

\section{The main construction}

\subsection{A particular Poisson manifold---Notations}

Let $(M, \omega)$ be a compact symplectic manifold. Let $\Omega\in
\mathcal{C}^\infty(]-\epsilon,\epsilon[,\Omega^2(M))$ be a smooth path of symplectic 
structures on $M$ such that $\Omega_{(0)}=\omega$. The smooth family $\{
\Omega{(t)}\}_{t\in]-\epsilon,\epsilon[}$ then canonically defines on 
$\hat{M}:=M\times]-\epsilon,\epsilon[$ a Poisson structure $\hat{\Omega}$ 
whose symplectic leaves are $\{(M\times\{t\},\Omega{(t)})\}$. We will 
denote by $\d\subset T\hat{M}$ the characteristic distribution of the 
Poisson structure $\hat{\Omega}$ (i.e. $\d_{(x,t)}=T_{(x,t)}(M\times\{t\})$).

\noindent The spaces $\mathcal{C}^\infty(\hat{M})[[\hbar]]$ (resp. 
$\mathcal{C}^\infty(M)[[\hbar]]$)
of power series in $\hbar$ with values in the algebra of smooth functions
on $\hat{M}$ (resp. $M$) are  
$\mathbb{R}[[\hbar]]$-algebras. The quotient 
$\mathbb{R}[[\hbar]]$-algebra $\mathcal{C}^\infty(M)[[\hbar]]/\hbar^{n+1}\mathcal{C}^\infty(M)[[\hbar]]$
will be denoted by $\mathcal{C}^\infty(M)[[\hbar]]_n$. 
It will often be identified with the space
of polynomials in $\hbar$ of degree at most $n$ with values in 
$\mathcal{C}^\infty(M)$.\\
We will consider the  natural inclusion 
$i:\mathcal{C}^\infty(M)[[\hbar]]_n \hookrightarrow 
\mathcal{C}^\infty(\hat{M})[[\hbar]]$ defined by $i(f)(x,t)=f(x) \ \forall
t \in ]-\epsilon , \epsilon[$. We will often denote $i(f)$ by 
$\hat{f}$.

\noindent By $DO_{\mathcal{D}}(\hat{M})$ we will denote the algebra of tangential 
(with respect 
to the distribution $\mathcal{D}$) differential 
operators on $\hat{M}$, i.e. the set of all differential operators on $\hat{M}$
vanishing on leafwise constant functions. By
$DO(M)$ we will denote the algebra of differential 
operators on $M$. As above we can consider the $\mathbb{R}[[\hbar]]$-algebras, 
$DO_{\mathcal{D}}(\hat{M})[[\hbar]]$ and
$DO(M)[[\hbar]]/\hbar^{n+1}DO(M)[[\hbar]]$ (abbreviated by $DO(M)[[\hbar]]_n$).
When dealing with bidifferential operators, we will use the prefix ``$biDO$".

\subsection{Taylor expansions}
We have 
$\mathcal{C}^\infty(\hat{M})
\simeq 
\mathcal{C}^\infty(]-\epsilon ,\epsilon[ ,\mathcal{C}^{\infty}(M))$ 
seeing every element $a\in \mathcal{C}^\infty(\hat{M})$ 
as a function of one variable with values
in a Fr\'echet space. We can therefore consider \cite{BkvarI} its Taylor 
expansion of order $n$ at $0$:
$$a(t)=\sum_{k=0}^{n}t^k\ \frac{1}{k!}\ a^{(k)}(0) + t^n R_n(u)(t) \, 
\hbox{with} \, R_n(u)(t)\rightarrow 0 \, \hbox{as} \, t\rightarrow 0.$$
We define the  $\mathbb{R}$-linear map, 
$j_n^{\hbar}:\mathcal{C}^\infty(\hat{M})\rightarrow
\mathcal{C}^\infty(M)[[\hbar]]_n$ by 
$j_n^{\hbar}a=\sum_{k=0}^{n}\hbar^k \frac{1}{k!}a^{(k)}(0)$.
It is extended to $\mathcal{C}^\infty(\hat{M})[[\hbar]]$ in the
following way:

\begin{eqnarray*}
  \mathcal{C}^\infty(\hat{M})[[\hbar]] & \longrightarrow 
& \mathcal{C}^\infty(M)[[\hbar]]_n\\
  a=\sum_{l\geqslant 0}\hbar^l a_l & \longmapsto 
& j_n^{\hbar}a
       =\sum_{k=0}^{n}\hbar^l j^{\hbar}_{n-k}a_l
       =\sum_{0\leqslant k+l\leqslant n}\hbar^{k+l}\frac{1}{k!} \ 
       a_l^{(k)}(0).
\end{eqnarray*}
One then has
\begin{lem} \label{hom} \hfill 
\begin{enumerate}
\item $j_n^{\hbar}a=j_n^{\hbar}(a\ \hbox{\tt mod}\ \hbar^{n+1})$.
\item $j_n^{\hbar}$ is an $\mathbb{R}[[\hbar]]$-algebra homomorphism.
\end{enumerate}
\end{lem}
\noindent We now extend the map $j_n^{\hbar}$ to 
$DO_{\mathcal{D}}(\hat{M})[[\hbar]]$
in the  natural way. 

\begin{dfn} \hfill
\begin{enumerate}
\item For $\Phi \in DO_{\mathcal{D}}(\hat{M})[[\hbar]]$, we define the 
operator 
$j_n^{\hbar}\Phi:\mathcal{C}^\infty (M)[[\hbar]]\to\mathcal{C}^\infty 
(M)[[\hbar]]_n$ 
by
$$j_n^{\hbar}\Phi\ .\ f = j_n^{\hbar}(\Phi .\hat{f})\ ,
\ \forall f\in \mathcal{C}^\infty (M)[[\hbar]].
$$
\item Similarly, for $B \in biDO_{\mathcal{D}}(\hat{M})[[\hbar]]$, we set 
$j_n^{\hbar} B\ .\ (f,g) = j_n^{\hbar}(B .(\hat{f}, \hat{g}))
\ ,\ \forall f,g \in \mathcal{C}^\infty (M)[[\hbar]]$.
\end{enumerate}
\end{dfn}
\begin{lem}\label{op} \hfill
\begin{enumerate}
\item
One has $j_n^{\hbar}\Phi \in DO(M)[[\hbar]]_n$ and $j_n^{\hbar} B \in biDO(M)[[\hbar]]_n$.
\item For all $a,b\in \mathcal{C}^\infty (\hat{M})[[\hbar]]$ one has 
$j_n^{\hbar}(\Phi .a)=j_n^{\hbar}\Phi\ .\ j_n^{\hbar}a$ and $j_n^{\hbar}(B.(a,b))=j_n^{\hbar}B\ .\ (j_n^{\hbar}a,j_n^{\hbar}b)$.
\end{enumerate}
\end{lem}
\Pf We will show that $j^\hbar_n\Phi$ and $j^\hbar_n B$
are local hence differential
 by Peetre's theorem \cite{Pe59,Pe60,CGDW80a}.
Let $f\in\mathcal{C}^\infty (M)$ and $U$ be an open set in $M$ 
such that $f_{/U}\equiv 0$. 
Since
$\hat{f}(x,t)~=~0\ \forall (x,t)\in U\times]-\epsilon , \epsilon[$ and 
 $\Phi$ is differential, one has
$(\Phi.\hat{f})_{/U\times]-\epsilon , \epsilon[}\equiv0$.
Hence
$$(j^\hbar_n\Phi) .f)_{/U}=(j^\hbar_n(\Phi .\hat{f}))_{/U}=
\sum_{0\leqslant k+l\leqslant n}\frac{\hbar^{k+l}}{l!}
(\Phi .\hat{f})^{(l)}_{/U\times]-\epsilon , \epsilon[}(0)=0.$$
The bidifferential case follows in the same way. This proves the first 
part of the lemma. The second one follows from simple computations.
\EPf
\begin{rmk}
{\rm 
Lemma \ref{op} implies that
$j_n^{\hbar}$, defined as a map from
$DO_{\mathcal{D}}(\hat{M})[[\hbar]]$ to $DO(M)[[\hbar]]_n$, is an
$\mathbb{R}[[\hbar]]$-algebra homomorphism for the composition product
on both algebras.
}
\end{rmk}
\subsection{Induced star-products}
Let now $\hat{\star}$ be any tangential star product on 
$(\hat{M},\hat{\Omega})$; for instance consider the Moyal-Fedosov star 
product defined in Section~\ref{FEDO}.
\begin{dfn}\label{def*} \hfill
\begin{enumerate}
\item We define $\star_n$ to be the operation from $\mathcal{C}^\infty (M)[[\hbar]]_n\times
\mathcal{C}^\infty (M)[[\hbar]]_n$ to $\mathcal{C}^\infty (M)[[\hbar]]_n$
given by 
$$f\star_n g=j_n^{\hbar}(\hat{f}\hat{\star}\hat{g}).$$ 
Equivalently
(by Lemma \ref{op}), seeing $\hat{\star}$ as an element of 
$biDO_{\mathcal{D}}(\hat{M})[[\hbar]]$, one has $\star_n=j_n^{\hbar}\hat{\star}.$
\item We define $\star$ to be  the operation from $\mathcal{C}^\infty (M)[[\hbar]]\times
\mathcal{C}^\infty (M)[[\hbar]]$ to $\mathcal{C}^\infty (M)[[\hbar]]$
given by
$$f\star g \ \hbox{\tt mod}\ \hbar^{n+1}
=f\ \hbox{\tt mod}\ \hbar^{n+1}\ \star_n\ g\ \hbox{\tt mod}\ \hbar^{n+1}$$
for all $n$ in $\mathbb{N}.$
\end{enumerate}
\end{dfn}
\begin{lem}\label{SP} \hfill
\begin{enumerate}
\item $\star_n$ is an associative product on the $\mathbb{R}[[\hbar]]$-algebra 
                                         $\mathcal{C}^\infty (M)[[\hbar]]_n.$
\item $\star$ is a star-product on $M$, called the {\bf induced star 
product} on $M$ by $\hat{\star}$.
\end{enumerate}
\end{lem}
\Pf 
For $f,g,h \in\mathcal{C}^\infty (M)[[\hbar]]_n$, one has
$(\hat{f} \hat{\star} \hat{g}) \hat{\star} \hat{h}=\hat{f} 
\hat{\star} (\hat{g} \hat{\star}\hat{h})$. Therefore,
$j_n^{\hbar}(\hat{f} \hat{\star} \hat{g}) \hat{\star} \hat{h}=j_n^{\hbar}\hat{f} \hat{\star} (\hat{g} \hat{\star}\hat{h})$
if and only if
\begin{eqnarray*}      
j_n^{\hbar}\left( 
        \hat{\star}.\left(\hat{\star}.(\hat{f},\hat{g}),\hat{h}\right) \right)
        =j_n^{\hbar}\left(
        \hat{\star}.\left(\hat{f},\hat{\star}.(\hat{g},\hat{h}) \right)\right)
        \mbox{ (reformulation)}  \\
\Leftrightarrow (j_n^{\hbar}\hat{\star}).\left(
        j_n^{\hbar}\left(\hat{\star}.(\hat{f},\hat{g})\right),j_n^{\hbar}\hat{h}
                              \right)                             
        =(j_n^{\hbar}\hat{\star}).\left(
        j_n^{\hbar}\hat{f},j_n^{\hbar}\left(\hat{\star}.(\hat{g},\hat{h})\right)
                                \right)
        \mbox{ (by Lemma \ref{op})}  \\
\Leftrightarrow 
       (j_n^{\hbar}\hat{\star}).\left((j_n^{\hbar}\hat{\star}).(f,g),h\right)
       =(j_n^{\hbar}\hat{\star}).\left(f,(j_n^{\hbar}\hat{\star}).(g,h)\right)
       \mbox{ (by Lemma \ref{op})}  \\
\Leftrightarrow (f\star_ng)\star_n h = f \star_n (g\star_n h)
       \mbox{ (by Definition \ref{def*}).} 
\end{eqnarray*}
This proves item 1 which is a classical way to show that a star-product is 
associative.
\EPf
\begin{cor}
If $\ \hat{\star}_1$ and $\hat{\star}_2$ are tangentially equivalent 
tangential star products on $(\hat{M},\hat{\Omega})$, then the induced star 
products $\star_1$ and $\star_2$ on $(M,\omega)$ are equivalent.
\end{cor}
\Pf The hypothesis implies that there exists an equivalence  
$\Phi \in DO_{\mathcal{D}}(\hat{M})[[\hbar]]$ such that
$\Phi.\left(a\hat{\star}_1b\right)=\Phi.a\ \hat{\star}_2\ \Phi.b$ for 
all $a,b \in \mathcal{C}^\infty (\hat{M})[[\hbar]]$. We then check, 
as in the proof of Lemma \ref{SP} that the operator 
$\Psi\ \hbox{\tt mod}\ \hbar^{n+1}:=j_n^{\hbar}\Phi\ ,\  n\in\mathbb{N}$ 
defines an equivalence between $\star_1$ and $\star_2$.
\EPf
\section{Characteristic classes}
Let $\ds\Omega^\hbar=\sum_{k\geqslant 0}\hbar^k\omega^k\in Z^2(M)[[\hbar]]$ 
be a formal power series of closed 2-forms on $M$. A refinement of
the classical Borel lemma (see the appendix) yields
\begin{lem}\label{borel}
Let $\ds\Omega_i^\hbar \in Z^2(M)[[\hbar]]\quad (i=1,2)$. Assume that 
$[\Omega^\hbar_1]=[\Omega_2^\hbar]$ in $H^2(M)[[\hbar]]$ or,
equivalently, that there exists $\nu^{\hbar}\in \Omega^1(M)[[\hbar]]$
such that $\Omega_2^\hbar - \Omega^\hbar_1 = d\nu^{\hbar}$.
Then there 
exists smooth functions $\Omega_i\in \mathcal{C}^\infty(]-\epsilon,\epsilon[,
\Omega^2(M))$ and $\nu\in \mathcal{C}^\infty(]-\epsilon,\epsilon[, \Omega^1(M))$ such that
\begin{enumerate}
\item[(i)] $\frac{1}{k!}\frac{d}{dt}\Omega_i|_{t=0}=\omega_i^k$;
\item[(ii)] $\forall t, \,\Omega_i(t)$ is symplectic;
\item[(iii)] $\forall t, \, \Omega_2(t) - \Omega_1(t) = d(\nu (t))$
or, equivalently, $[\Omega_1(t)]=[\Omega_2(t)]$.
\end{enumerate}
\end{lem}
\begin{dfn}
Let us fix a connection $\nabla^{0}$ in the vector bundle 
$$\d\to\hat{M}=M\times]-\epsilon,\epsilon[.$$
Let $\Omega^\hbar\in\Omega^2(M)[[\hbar]]$ be a series of closed 2-forms 
on $M$ such that $\Omega^\hbar\mbox{ mod }\hbar=\omega$. Let $\Omega\in 
\mathcal{C}^\infty(]-\epsilon,\epsilon[,\Omega^2(M))$ be a smooth family of 
symplectic structures on $M$ admitting $\Omega^\hbar$ as $\infty$-jet (cf. 
Lemma~\ref{borel}). Let $\nabla$ be the symplectic foliated connection on $\hat{M}$ 
obtained from the data of $\nabla^0$ and $\Omega$ (cf. Section~\ref{FEDO}). Let 
$\hat{\star}$ be the Moyal-Fedosov star product on 
$(\hat{M},\hat{\Omega})$ associated to $\nabla$. The 
star product $\star^{\Omega^\hbar}$ on $(M,\omega)$ induced by 
$\hat{\star}$ will be called {\bf the star 
product associated to the series $\Omega^\hbar$}.
\end{dfn}
\begin{prop}\label{CC}
Let $\Omega^\hbar_i\,(i=1,2)$ be two series of closed 2-forms on $M$ such 
that $\Omega^\hbar_i\mbox{ mod }\hbar=\omega.$ Denote by $\star_i\,(i=1,2)$ 
the associated star products on $(M,\omega)$. Then $\star_1$ and $\star_2$ 
are equivalent star products if and only if 
$[\Omega^\hbar_1]=[\Omega^\hbar_2]$ 
in $H_{\mbox{de Rham}}^2[[\hbar]]$.
\end{prop}
\noindent The proof of Proposition~\ref{CC} is postponed to the end of
this section.
\begin{dfn}\label{FF}
A diffeomorphism $\hat{\varphi}~: \hat{M}\to\hat{M}$ {\bf preserves the foliation}
if 
\begin{enumerate}
\item[(i)] $\hat{\varphi}(M_t)\subset M_t\quad \forall t$ and
\item[(ii)] $\hat{\varphi}|_{M_0}=id_{M_0}$.
\end{enumerate}
\end{dfn}
\noindent We first adapt Moser's lemma to our parametric situation.
\begin{lem}
Let $\{\Omega_i(t)\}_{t\in]-\epsilon,\epsilon[}\,(i=1,2)$ be two smooth 
families of symplectic structures on $M$ such that 
$\Omega_1(0)=\Omega_2(0)=\omega$. Assume that, for all 
$t\in]-\epsilon,\epsilon[$ they have the same de Rham class~:
$
[\Omega_1(t)]=[\Omega_2(t)]
$
in $H^2(M)$. Then there exists a Poisson diffeomorphism $\hat{\varphi}~: 
(\hat{M},\hat{\Omega_2})\to
(\hat{M},\hat{\Omega_1})$ which preserves the foliation.
\end{lem}
\Pf By Hodge's theory one has that 
$\Omega_{1}(t)-\Omega_{2}(t)=d\nu^t$ where $ \nu^t\in \Omega^1(M)$ is smooth in $t$.  
Set $\omega_s^t=\Omega_2(t)+s\,d\nu^t\ , \ s\in[0,1]$. The form $\omega_s^0=\omega$ 
is symplectic on $M$ for all $s\in[0,1]$; hence by compactness, one can 
choose $\epsilon>0$ such that $\omega_s^t$ is symplectic for all 
$t\in]-\epsilon,\epsilon[$ and $s\in[0,1]$.\\
Consider $N=M\times[0,1]$ endowed with the natural 
foliation $\mathcal{F}=\{M\times\{s\}\}$.
Define the following smooth families of 2-forms on $N$~:
\begin{eqnarray*}
(\tilde{\omega}_t)_{(x,s)}:=(\omega^t_s)_x\mbox{ and}\\
(\omega_t)_{(x,s)}:=(\tilde{\omega}_t)_{(x,s)}-(\nu^t)_x\wedge ds.
\end{eqnarray*}
Then $d_N(\omega_t)=d_N(\tilde{\omega})-d_M(\nu^t)\wedge ds=0$ for all $t$. 
Moreover,
$\mbox{rad}_{T_{(x,s)}(N)}(\omega_t)$ is not entirely contained in 
$T(\mathcal{F})$;
hence one can find a smooth family of vector fields of the form~: 
$X_t=\frac{\partial}{\partial_s}+Y_t \quad (Y_t\in\Gamma(T(\mathcal{F})))$ generating the smooth 
family of smooth distributions~: $\mbox{rad}(\omega_t)$.\\
One has therefore 
$$
\mathcal{L}_{X_t}\omega_t=d( i_{X_t}\omega_t ) +i_{X_t}d\omega_t=0.
$$
Denoting by $\{\varphi_{X_t}^u\}$ the flow of $X_t$, one has~:
\begin{eqnarray*}
( \varphi_{X_t}^u )^\star\omega_t=\omega_t\mbox{ and}\\
\varphi_{X_t}^u(M\times\{s\})=M\times\{s+u\}.
\end{eqnarray*}
One then gets a smooth family $\{\varphi_t\}$ of diffeomorphisms of $M$  
defined by
$$ 
\varphi_{X_t}^1\circ i_0=i_1\circ\varphi_t
$$
such that $\varphi_t^\star(\Omega_1(t))=\Omega_2(t)$ ($i_s:M\to N$ denotes 
the natural inclusion $i_s(x)=(x,s)$).\\
Shrinking $\epsilon$ once more if necessary, one gets the desired Poisson 
map by setting $\hat{\varphi}(x,t)=(\varphi_t(x),t)$. Observe that 
$X_0=\partial_s$, hence $\varphi_0=id_M$.
\EPf

\begin{lem}
Let $\hat{\star}_i\quad(i=1,2)$ be tangential star products on $\hat{M}$. 
Suppose there exists a diffeomorphism $\hat{\varphi}~: \hat{M}\to\hat{M}$ 
preserving the foliation such that 
$\hat{\star}_1^{\hat{\varphi}}=\hat{\star}_2\mbox{ mod }(\hbar^n)$. Then, 
$\star_1$ and $\star_2$ are equivalent star products up to order~$n$.
\end{lem}
\Pf
The right action of the diffeomorphism group, $\mathcal{C}^\infty(\hat{M})\times
\mbox{Diff}(\hat{M})\stackrel{\rho}{\to}\mathcal{C}^\infty(\hat{M})$,
$\rho(\hat{\varphi})u={\hat{\varphi}}^\star u$ yields a map~:
$$\rho_n^\hbar:\mbox{Diff}(\hat{M})\to 
\mbox{Hom}_{\R}(\mathcal{C}^\infty(M),\mathcal{C}^\infty(M)[[\hbar]]_n)~:
\rho_n^\hbar(\hat{\varphi})f=\mbox{j}^\hbar_n(\hat{\varphi}^\star\hat{f}).$$
Definition~\ref{FF} implies that if $\hat{\varphi}$ preserves the 
foliation, then $\rho_n^\hbar(\hat{\varphi})\in \mbox{DO}(M)[[\hbar]]_n$ and 
$\rho_0^\hbar(\hat{\varphi})=id$. Therefore an argument similar to the one used 
for Lemma~\ref{SP} yields the conclusion.
\EPf

\begin{cor}
Within the notations of Proposition~\ref{CC}, if $\Omega_1^\hbar$ and
$\Omega_2^\hbar$ are cohomologous in $H^2(M)[[\hbar]]$, 
then the star products $\star_1$ and $\star_2$ are equivalent.
\end{cor}
\Pf
The first $N$ cochains of a Fedosov star product are entirely determined by the 
$N$ first terms of its Weyl curvature. Therefore, the above lemmas imply that 
$\star_1$ and $\star_2$ are equivalent up to any order. It is then classical 
that they are equivalent \cite{BFFLS77a}.
\EPf
\noindent {\sl Proof of Proposition~\ref{CC}}

\noindent We first consider a particular case. Let  
$\alpha^{\hbar}=\alpha^{\rm o}+\hbar\alpha^1\ldots \in
Z^2(M)[[\hbar]]$ be a sequence of closed 2-forms on $M$. Set
$\Omega^{\prime\hbar} =\Omega^{\hbar} + \hbar^k\alpha^{\hbar}$.
Denote by $\Omega,\,\alpha$ and $\Omega'=\Omega+t^k\alpha$
respectively the smooth functions associated to the series
$\Omega^\hbar, \alpha^\hbar$ and $\Omega^{\prime\hbar}$ as in
Lemma~\ref{borel}. The function $\Omega'$ defines a Poisson structure
on $\hat{M}$. We denote by $\Lambda'$ (resp. $\omega'$) the
corresponding bivector field (resp. $\d$-2-form). One has
\begin{equation}\label{dif}
\omega^{\prime t} =\omega^{t} + t^k\alpha^{t}\quad \hbox{and}\quad
\Lambda^{\prime t} = \Lambda^t - t^k \sharp\alpha^{\rm o} +t^{k+1} \lambda\ ,
\end{equation}
where we denote again by $\alpha^{t}$ the $\mathcal{D}$-2-form corresponding 
to $\alpha^{t}$ and where
$\lambda$ is an element of $\mathcal{C}^{\infty}(]-\epsilon ,
\epsilon[,
\Gamma\wedge^2\mathcal{D})$.
Let $\star^\prime$ be the star-product on $M$ induced by the 
Moyal-Fedosov star-product $\hat\star^\prime$ on $(\hat{M}, \Lambda^{\prime})$.
We now define a specific
foliated symplectic connection $\nabla^\prime$ adapted to $\omega^{\prime t}$.
Let us look
for $\nabla^\prime$ of the form $\nabla + S$ where $S$ is a symmetric
2-$\mathcal{D}$-tensor field. We set
$$
\omega^{\prime t}\left(\nabla^\prime_u v , w \right)
= \omega^{\prime t}\left(\nabla_ u v , w \right)
+\frac{1}{3}\left(\nabla_u\omega^{\prime t}\right)(v,w)
+\frac{1}{3}\left(\nabla_v\omega^{\prime t}\right)(u,w). 
$$
This leads to
$\left(\omega^t + t^k\alpha\right)\left(S(u,v), w \right)
=\frac{t^k}{3}\left[\left(\nabla_u\alpha\right)(v,w)
+\left(\nabla_v\alpha\right)(u,w)\right]$ as $\nabla\omega^t=0$.
By construction $\omega^t + t^k\alpha^t$ is invertible, so $S(u,v)$ is completely
determined and of the form $S(u,v)=t^ks(u,v)$. We thus have
\begin{equation}\label{conn}
\nabla^\prime = \nabla + t^k s.
\end{equation}
Let now $\circ^t$ (resp. $\circ^{\prime t}$) be the associative
product on the sections of the Weyl bundle $\w$ over $\hat{M}$
determined by the data of $\Lambda$ (resp. $\Lambda'$)
(cf. Section~\ref{FEDO} and Remark~\ref{WB}). By construction, we then get
$\forall u, v \in \mathcal{W}$,
\begin{equation} \label{moyal}
\frac{d^l}{dt^l}(u\circ^t v - u\circ^{\prime t} v)(0)=0\quad \forall l\leqslant k-1.
\end{equation}
Similarly for Moyal-Fedosov star products, $\hat{\star}$ and
$\hat{\star}^{\prime}$, associated to $(\Omega,\nabla)$ and $(\Omega',
\nabla')$, (\ref{conn}) and (\ref{moyal}) yield
\begin{equation}\label{nul-k}
\frac{d^l}{dt^l}(a\hat{\star} b - a\hat{\star}^{\prime} b)(0)=0\quad \forall l\leqslant k-1.
\end{equation}
Now let us see what happens for $\star$ and $\star^\prime$. Let  $f, g \in \mathcal{C}^\infty(M)$
and write $\hat{\star}=\sum_{i\geqslant 0}\hbar^i \hat{C}_i$ 
and $\hat{\star}^{\prime}=\sum_{i\geqslant 0}\hbar^i \hat{C}_i^\prime$.
Setting $u^{(l)}:=\frac{d^l}{dt^l}u$, we have
\begin{eqnarray*}
& & f\star g - f\star^\prime g = \\
& = & \sum_{j\geqslant 0} 
\frac{\hbar^j}{j!}(f\hat{\star} g - f\hat{\star}^\prime g)^{(j)}(0)
 =  \sum_{j\geqslant k} \frac{\hbar^j}{j!}(f\hat{\star} g - f\hat{\star}^\prime g)^{(j)}(0)
\qquad \mbox{ (cf. equation (\ref{nul-k})) }  \\
& = & \sum_{j\geqslant k} \frac{\hbar^j}{j!} \left( \sum_{i\geqslant 0} \hbar^i
\left(  \hat{C}_i(f,g) - \hat{C}_i^\prime (f,g)\right) \right)^{(j)}(0)\\
& = & \sum_{j\geqslant k, i\geqslant 0} \frac{\hbar^{i+j}}{j!}
\left( \hat{C}_i(f,g) - \hat{C}_i^\prime (f,g) \right)^{(j)}(0)\\
& = & \sum_{m\geqslant k} \hbar^m 
\sum_{m=i+j, j\geqslant k, i\geqslant 0}\frac{1}{j!} \left( \sum_{k\geqslant 0} 
\left( \hat{C}_i(f,g) - \hat{C}_i^\prime (f,g)\right) \right)^{(j)}(0) \\
& = & \frac{\hbar^k}{k!} \left( fg - gf \right)
      + \frac{\hbar^{k+1}}{(k+1)!} \left( fg - gf \right)
      + \frac{\hbar^{k+1}}{k!} \left( \hat{C}_1(f,g) - \hat{C}_1^\prime (f,g) \right)\\
& & + \sum_{m\geqslant k+2} \hbar^m 
        \sum_{m=i+j, j\geqslant k, i\geqslant 0}\frac{1}{j!} \left( \sum_{k\geqslant 0} 
        \hat{C}_i(f,g) - \hat{C}_i^\prime (f,g) \right)^{(j)}(0) \\
& = & \frac{\hbar^{k+1}}{k!} \left( \Lambda^t(df,dg) - \Lambda^{\prime t}(df,dg) \right)^{(k)}(0)
      + o(\hbar^{k+1})\\
& = & \frac{\hbar^{k+1}}{k!} \left( \Lambda^t(df,dg) 
        - \Lambda^t (df,dg) + t^k \sharp\alpha^{\rm o}(df,dg) -
     t^{k+1} \lambda (df,dg) \right)^{(k)}(0) \\
&  &  + o(\hbar^{k+1})\\
& = & \hbar^{k+1} \sharp\alpha^{\rm o}(df,dg) + o(\hbar^{k+1})
\end{eqnarray*}

\noindent Then, setting $\star=\sum_{i\geqslant 0}\hbar^i C_i$ 
and $\star^\prime=\sum_{i\geqslant 0}\hbar^i C_i^\prime$, we have 
\begin{equation}
C_i^\prime = C_i \; ,\quad i=0,\ldots , k \quad \mbox{and} \quad
C_{k+1}^\prime =  C_{k+1} + \sharp\alpha^{\rm o}.
\end{equation}
Let us pass to the general case. Suppose that 
$[\Omega^{\hbar}_1]\neq[\Omega^{\hbar}_2]$. We denote by $k$ the smallest
integer such that $[\omega_1^k]\neq[\omega_2^k]$.
Let us consider $\Omega^{\hbar}_3=\hbar\omega^1_1 + \hbar^2 \omega^2_1 +\cdots 
+ \hbar^{k-1}\omega_1^{k-1}+\hbar^m\omega_2^{k}+\hbar^{k+1}\omega_2^{k+1}+\cdots$.
We have $[\Omega^{\hbar}_3]=[\Omega^{\hbar}_2]$ and 
$\Omega^{\hbar}_1=\Omega^{\hbar}_3 +\hbar^k(\omega_1^{k} - \omega_2^{k})
+\hbar^{k+1}\cdots$.
Denoting by $\star_{i}$ the product associated with 
$=\Omega^{\hbar}_i$, we know that $\star_2$ and $\star_3$ are equivalent.
What has been done previously implies $\star_1=\star_3\ \hbox{\tt mod}\ \hbar^{k+1}$
and $C^{(1)}_{k+1} =  C^{(3)}_{k+1} \pm \sharp\alpha^{\rm o}$ with
$\alpha^{\rm o}=\omega_1^{k} - \omega_2^{k}$. But in this case, we know
that $\star_1 \sim \star_3\ \hbox{\tt mod}\ \hbar^{k+2} $ if and only if
$\alpha^{\rm o}$ is exact \cite{BCG97a}. Since  $\omega_1^{k} - \omega_2^{k}$
is not exact by hypothesis, $\star_1 \not\sim \star_3$ and thus
$\star_1 \not\sim \star_2$.
\EPf
\section{Appendix: Borel's Lemma}

\begin{prop}[Borel's Lemma] \label{borAP}\hfill \\
Let $M$ be a compact smooth manifold of dimension $d$ and 
$\{\alpha_n\in\Omega^q(M)\ |\ n\in\mathbb{N}\}$ be a sequence of 
q-forms on $M$.
Then there exists $f \in\mathcal{C}^\infty (\mathbb{R},\Omega^q(M))$
such that $\ds\frac{d^nf}{dt^n}(0)=\alpha_n$.
\end{prop}

\medskip

\Pf
Let $\varphi\in\mathcal{C}^\infty (\mathbb{R},\mathbb{R})$ be
such that $\varphi (t)=1$ for $|t|\leqslant\frac{1}{2}$ and 
$\varphi (t)=0$ for $|t|\geqslant 1$ and set
$f_n:\mathbb{R}\longrightarrow\Omega^q(M)\quad \ds f_n(t)=\frac{\alpha_n}{n!}\varphi(\lambda_nt)$
where the numbers $\{\lambda_n\}$ will be defined later.
Let $\ds\{V_i, \psi_i, i=1,\ldots , N\}$ be a finite ($M$ is compact)
atlas of $M$ trivializing the bundle $\wedge^qT^*M\rightarrow M$.
Restricted to such a chart, we can view a section of $\wedge^qT^*M\rightarrow M$
as a function from $\psi_i(V_i)\rightarrow \mathbb{R}^{s}$ with
$s=\cbi{d}{q}$. We denote
by $|| . ||$ the Euclidean norm on $\mathbb{R}^{s}$. One can make a 
choice of $\{\lambda_{n}\}$ is such a way that, for all
$\ds k\in \mathbb{N}$ such that $0\leqslant k\leqslant n-1$, one has
$\sup\{||D^\nu . \frac{\partial^l}{\partial t^l} . f_n(x,t)||\ 
|\ (x,t)\in V_i\times\mathbb{R},\ |\nu |+l\leqslant k,\ i=1,\ldots , N\}
\leqslant\frac{1}{2^n}$
where $\nu =(\nu_1 , \ldots ,\nu_d)\in\mathbb{N}^{d}$,
$|\nu|=\nu_1 +\ldots +\nu_d$ and 
$D^\nu=\ds\frac{\partial^{\nu_1}}{\partial x_1^{\nu_1}}
\ldots\frac{\partial^{\nu_d}}{\partial x_d^{\nu_d}}$. Indeed, let us fix a $V_i$. In this chart we have 
$$\ds D^\nu . \frac{\partial^l}{\partial t^l} . f_n(x,t)
=\sum_{p=0}^{l}\cbi{p}{l}\frac{n!}{(n-p)!}\frac{D^\nu\alpha_n}{n!}t^{n-p}\lambda_n^{l-p}
\varphi^{(l-p)}(\lambda_nt).$$
Define $\ds K_n=\sum_{|\nu|=1}^{n}\sup_{x\in V_i}||D^\nu\alpha_n(x)||$
and $\ds M_n=\sum_{j=1}^{n}\sup_{t\in \mathbb{R}}|\varphi^{(j)}(t)|$.
On the support of $\varphi^{(l-p)}(\lambda_nt)$ we have $\lambda_nt\leqslant1$.
Hence $\ds||D^\nu . \frac{\partial^l}{\partial t^l} . f_n(x,t)||
\leqslant\sum_{p=0}^{l}\cbi{p}{l} \frac{K_n}{(n-p)!}\frac{M_n}{\lambda_ n^{n-l}}
\leqslant\frac{K_n M_n}{\lambda_ n}\sum_{p=0}^{n-1}\cbi{p}{n-1}\frac{1}{(n-p)!}$
as $n-p\geqslant n-l\geqslant n-k\geqslant 1$. Therefore $\lambda_ n^{n-p}\geqslant\lambda_n$ if 
$\lambda_n\geqslant 1$ and $\cbi{p}{l}\leqslant \cbi{p}{n-1}$.
Thus a choice of the $\lambda_{n}$'s such that  
$$\ds\lambda_n\geqslant\max\{1,2^nK_nM_n\sum_{p=0}^{n-1}\cbi{p}{n-1}\frac{1}{(n-p)!}\}$$
yields the assertion on $V_i$. The conclusion follows by finiteness 
of our atlas. In particular the function $\ds 
f:=\sum_{n=0}^{\infty}f_n$
is well defined. By the preceding lemma $f$ is well defined. Moreover,
$f\in\mathcal{C}^k(M\times\mathbb{R},\wedge^qT^*M)$ for all $k\in 
\mathbb{N}$ as it appears when writing $\ds 
f(x,t)=\sum_{n=0}^{k}f_n(x,t) + \sum_{n=k+1}^{+\infty}f_n(x,t)$. 
Therefore $f\in \mathcal{C}^\infty (M\times\mathbb{R},\wedge^qT^*M)$.
By definition of $f$, $f(x,t)\in\wedge^qT_x^*M$, hence $f( . ,t)$ is a smooth
section of $\wedge^qT^*M\rightarrow M$ and 
$f:\mathbb{R}\rightarrow\Omega^q(M)$ with $f(t)=f( . ,t)$ belongs
to $\mathcal{C}^\infty (\mathbb{R},\Omega^q(M))$. Moreover, we have 
$$\ds f^{(k)}(t)=\sum_{n=0}^k\frac{\alpha_n}{n!}\big(t^n\varphi(\lambda_nt)\big)^{(k)}
+ \sum_{n=k+1}^\infty\sum_{p=0}^k\cbi{p}{k}\frac{n!}{(n-p)!}
\frac{\alpha_n}{n!}t^{n-p}\varphi^{(k-p)}(\lambda_nt).$$
In the second sum, we have $n-p\geqslant n-k\geqslant 1$. Thus it vanishes for $t=0$.
In the first sum, if $n\leqslant k-1$, $\varphi$ is differentiated at least once. 
As $\varphi^{(j)}(0)=0$ for $j\geqslant 1$, it vanishes for $t=0$.
Therefore, we have
$\ds f^{(k)}(0)=\frac{\alpha_k}{k!}\big(t^k\varphi(\lambda_kt)\big)^{(k)}(0)
=\frac{\alpha_k}{k!}\sum_{p=0}^k\cbi{p}{k}\frac{k!}{(k-p)!}t^{k-p}\lambda_k^{k-p}
\varphi^{(k-p)}(\lambda_kt)|_{t=0}$.\\
For $p\leqslant k-1$, we have $k-p\geqslant 1$ and the corresponding term vanishes.
Hence  $\ds f^{(k)}(0)=\frac{\alpha_k}{k!}k!\varphi(0)=\alpha_k$.
\EPf

\begin{cor}\label{corborel} Let $(\alpha_n)_{n\in\mathbb{N}}$, 
$(\alpha^1_n)_{n\in\mathbb{N}}$, $(\alpha^2_n)_{n\in\mathbb{N}}$ and
$(\nu_n)_{n\in\mathbb{N}}$ be sequences of forms on $M$. Then there 
exist smooth functions $f^1$, $f^2$ and $f$ corresponding respectively to $(\alpha^1_n)_{n\in\mathbb{N}}$, 
$(\alpha^2_n)_{n\in\mathbb{N}}$ and $(\nu_n)_{n\in\mathbb{N}}$ as in Proposition 
\ref{borAP} such that
\begin{enumerate}
\item if\ $d\alpha_n=0,\ \forall n\in \mathbb{N}$ then,
          $d(f(t))=0,\ \forall t\in\mathbb{R}$.
\item if\ $\alpha^2_n - \alpha^1_n = d\nu_n\ \forall n\in\mathbb{N}$, then
          $f^2(t) - f^1(t) = d(f(t),\ \forall t\in\mathbb{R}$.
\end{enumerate}
\end{cor}

\Pf\hfill
1) We have $\ds f_n(t)=\frac{\alpha_n}{n!}\varphi(\lambda_nt)$ hence $d(f_n(t))=0,\ \forall t \in \mathbb{R}$.
For each $t$, $f(t)=\sum_{n=0}^{\infty}f_n(t)$ converges
absolutely in $\Gamma^1(M,\wedge^2T^*M)$.

2) Let $\lambda^1_n$, $\lambda^2_n$ and $\lambda_n$ be three real 
sequences defining smooth functions $\tilde{f}^1$, $\tilde{f}^2$ and $\tilde{f}$ corresponding
respectively to $(\alpha^1_n)_{n\in\mathbb{N}}$, $(\alpha^2_n)_{n\in\mathbb{N}}$ and
$(\nu_n)_{n\in\mathbb{N}}$ as in the proof of Proposition \ref{borAP}.
Consider the sequence $\mu_n = max\{\lambda^1_n , \lambda^2_n , \lambda_n \}$.
Replacing $\lambda^1_n$, $\lambda^2_n$ and $\lambda^1_n$ by $\mu_n$ we 
get new functions $f^1$, $f^2$ and $f$ again corresponding
respectively to $(\alpha^1_n)_{n\in\mathbb{N}}$, $(\alpha^2_n)_{n\in\mathbb{N}}$ and
$(\nu_n)_{n\in\mathbb{N}}$ such that
$f^2_n - f^1_n = df_n\ \forall n\in\mathbb{N}$.  Since for $t$ fixed the series
converge absolutely for the $\mathcal{C}^0$ norm on the forms, we obtain the result. 
\EPf

\newcommand{\etalchar}[1]{$^{#1}$}

{\sc 
Universit\'e Libre de Bruxelles, Brussels, Belgium}\\
{\sc email:} pbiel@ulb.ac.be\\
{\sc Universit\'e de Bourgogne, Dijon, France}\\
{\sc email:} bonneau@u-bourgogne.fr


\begin{thebibliography}{CGDW80}

\bibitem[BCG97]{BCG97a}
M.~Bertelson, M.~Cahen, and S.~Gutt.
\newblock Equivalence of star products.
\newblock {\em Class. Quantum Grav.}, 14(1A):A93--A107, 1997.

\bibitem[BFF{\etalchar{+}}77]{BFFLS77a}
F.~Bayen, M.~Flato, C.~Fronsdal, A.~Lichnerowicz, and D.~Sternheimer.
\newblock Quantum mechanics as a deformation of classical mechanics.
\newblock {\em Lett. Math. Phys.}, 1(6):521--530, 1977.

\bibitem[Bou67]{BkvarI}
Nicolas Bourbaki.
\newblock {\em \'{E}l\'ements de math\'ematique. {F}asc. {X}{X}{X}{I}{I}{I}.
  {V}ari\'et\'es diff\'erentielles et analytiques. {F}ascicule de r\'esultats
  ({P}aragraphes 1 \`a 7)}.
\newblock Hermann, Paris, 1967.
\newblock Actualit\'es Scientifiques et Industrielles, No. 1333.

\bibitem[CGDW80]{CGDW80a}
M.~Cahen, S.~Gutt, and M.~De~Wilde.
\newblock Local cohomology of the algebra of ${C}\sp{\infty }$\ functions on a
  connected manifold.
\newblock {\em Lett. Math. Phys.}, 4(3):157--167, 1980.

\bibitem[Fed94]{Fedb94a}
Boris~V. Fedosov.
\newblock A simple geometrical construction of deformation quantization.
\newblock {\em J. Differential Geom.}, 40(2):213--238, 1994.

\bibitem[Fed96]{Fedbook}
Boris~V. Fedosov.
\newblock {\em Deformation quantization and index theory}.
\newblock Akademie Verlag, Berlin, 1996.

\bibitem[Kon97]{Konm97b}
Maxim Kontsevich.
\newblock Deformation quantization of {P}oisson manifolds {I}, 1997.
\newblock {\tt q-alg/9709040}.

\bibitem[Pe59]{Pe59}
Jaak Peetre.
\newblock {Une caract\'erisation abstraite des op\'erateurs diff\'erentiels.}
\newblock {\em Math. Scand.}, 7:211--218, 1959.

\bibitem[Pe60]{Pe60}
Jaak Peetre.
\newblock {Rectification \`a l'article ``Une caract\'erisation abstraite des
  op\'erateurs diff\'erentiels''.}
\newblock {\em Math. Scand.}, 8:116--120, 1960.

\bibitem[Vai94]{Vais94a}
Izu Vaisman.
\newblock {\em Lectures on the geometry of {P}oisson manifolds}.
\newblock Birkh\"auser Verlag, Basel, 1994.

\end{thebibliography}
\end{document}